\input amstex

\magnification=1200

\centerline{\bf Volume and Angle Structures on 3-Manifolds}
\medskip
\centerline{\bf Feng Luo}

\medskip
\centerline{\it Abstract}

We  propose an approach to find constant curvature metrics on triangulated closed 3-manifolds using a finite dimensional variational method whose energy function is the volume. The concept of an 
angle structure on a tetrahedron and on a triangulated closed 3-manifold is introduced following  
the work of Casson, Murakami and Rivin. It is proved by A. Kitaev and the author that any closed 3-manifold has a triangulation supporting an angle structure.
The moduli space of all angle structures  on a triangulated 3-manifold
 is a bounded open convex polytope in a Euclidean space.
The volume of an angle structure is defined.  
Both the angle structure and the volume are natural generalizations of tetrahedra in the constant
 sectional curvature spaces and their volume.  It is shown that the volume functional can be extended continuously
to the compact closure of the moduli space.  In particular, the maximum point of the volume functional 
always exists in the compactification.
The main result shows that for a 1-vertex triangulation of a closed
 3-manifold if the volume function on the moduli space has a local maximum point,  then either the manifold admits a constant curvature Riemannian metric or the manifold contains a non-separating 2-sphere or real projective plane.

\bigskip
\noindent
{\bf \S 1. Introduction}
\medskip
\noindent
1.1.
The purpose of the announcement is to propose an approach to find constant curvature metrics on triangulated closed 3-manifolds using a
finite dimensional variational method whose energy function is the volume. It is a generalization of the corresponding
program introduced by Casson [La] and Rivin [Ri1] for finding hyperbolic metrics on 
compact ideally triangulated 3-manifolds with torus boundary.  It is also motivated by the
work of Murakami [Mu]. Very recently, Casson and Rivin's approach
was successfully carried out by Francois Gueritaud [Gu] to give a new proof of the existence of
hyperbolic structures on 1-holed torus bundle over the circle with Anosov monodromy.

By an \it angle structure \rm on a 3-simplex we mean an assignment of a number, called the \it dihedral angle\rm,   
to each edge of the 3-simplex so that dihedral angles at three edges sharing a common vertex are the inner angles of a spherical
triangle.  (Note that the similar concept introduced by Casson and Rivin requires that the vertex triangle be a Euclidean triangle).
Since a spherical triangle is determined by its inner angles subject to four linear inequalities, the moduli space
of all  angle structures on a 3-simplex, denoted by $  AS(3)$,  is an open bounded  convex polytope in $\bold R^6$. 
Examples of  angle structures on a 3-simplex are classical geometric tetrahedra, i.e., 
Euclidean, hyperbolic and spherical tetrahedra measured by  dihedral angles. However, not every angle structure is of this form.
We define the \it generalized volume \rm (or \it volume \rm for short) of an angle structure on a 3-simplex  
by generalizing the Schlaefli formula.
To be more precise, the Schlaefli formula for volume of spherical or hyperbolic tetrahedra says that  volume  can be defined
 by integrating
the Schlaefli 1-form which depends on the dihedral angles and the edge lengths. Furthermore, one finds the edge length
from the dihedral angles by using the Cosine Law twice.  We follow this path to define the volume of  an angle structure
by first defining
the edge lengths using the Cosine Law and then verifying the resulting Schlaefli 1-form is closed in the moduli space  $  AS(3)$.
The generalized volume is  the integration of the Schlaefli 1-form. 
In the case of classical geometric tetrahedra, the generalized volume coincides with classical geometric volume for spherical and 
hyperbolic tetrahedra and is zero for Euclidean tetrahedra.

By an \it angle structure \rm on a triangulated closed 3-manifold  $(M, T)$ we mean an assignment of
a number, called the dihedral angle, to each edge of each 3-simplex  in the triangulation $T$, so that (1) the assignment is an angle structure on each 3-simplex in $T$, (2) 
the sum of dihedral angles at each edge in $T$ is $2\pi$. The basic examples of angle structures are totally geodesic triangulations in
a 3-manifold with a constant curvature metric. The \it volume \rm of an  angle structure
on a triangulated manifold is defined to be the sum of the volume of
its 3-simplexes.  From the definition, it is clear that the moduli space of all angle structures on a fixed triangulated 3-manifold $(M, T)$,
denoted by $  AS(M, T)$, forms a bounded convex polytope in a Euclidean space.  
The main theorem  is the following.
\medskip

\noindent
{\bf Theorem 1.1.} \it Suppose $(M, T)$ is a  closed 3-manifold with a triangulation $T$ so that the moduli space of all angle
structures $AS(M, T)$ is non-empty. Then volume function can be extended continuously to the compact closure of the space $AS(M, T)$. If the 
 volume functional has a local maximum point in $AS(M, T)$,
 then, 

(a)  the manifold $M$ supports a constant sectional curvature Riemannian metric, or
 
(b) there is a normal surface of positive Euler characteristic in the triangulation $T$ so that
it intersects each 3-simplex in at most one normal disk.

In particular, if the triangulation $T$ has only one vertex, then the normal surface in case (b) is non-separating.
\rm
\medskip

\medskip

The existence of angle structures on a triangulation is a linear programming problem.  
\medskip
\noindent
{\bf  Problem 1.2.} \it Does every 
closed irreducible non-Haken non Seifert-fibered 3-manifold have a 1-vertex triangulation supporting an angle structure?
\rm
\medskip
We expect the  problem has an affirmative solution. See [JR] for more information on efficient triangulaions.
In [LT], a relationship between the existence of angle structure and the normal surface theory is established.

If we do not assume 1-vertex condition, then the following has been proved by A. Kitaev and myself.

\medskip
\noindent
{\bf Theorem 1.3.} \it For any closed 3-manifold $M^3$, there is a triangulation of $M^3$ supporting an angle
structure. \rm

\rm

By theorem 1.1, the maximal point of the volume function on the compact closure of $AS(M, T)$ always exists. If the maximum point is in $AS(M,T)$, then we conclude either the
manifold $M$ is geometric, or the triangulation admits a special normal surface of positive Euler
characteristic.  It is expected ([Lu3]) that the maximum point in the boundary $\partial
AS(M, T)$ will give rise either a geometric structure on the manifold, or a 
special normal surface of non-negative Euler characteristic in the triangulation. 

\medskip
\noindent
1.2. Using dihedral angles instead of edge lengths to parameterize a classical geometric tetrahedra seems
to have some advantages. 
 First of all, the Schlaefli formula suggests that dihedral angles are natural variables with respect to the volume consideration. Second, dihedral angle parameterization puts
all hyperbolic, spherical and Euclidean tetrahedra in one frame work. Third,
volume considered as a function of
dihedral angles can be extended continuously to the degenerated tetrahedra. This was conjectured
by John Milnor [Mi] and was established recently in [Lu1] and [GL] ([Ri2] has a new proof together with some generalizations). On the
other hand, volume considered as a function of edge lengths cannot be extended continuously
to degenerated tetrahedra. The use of dihedral angle parameterization has also appeared in the work of Murakami [Mu].

\medskip
\noindent
1.3.
In this subsection, we briefly sketch the main ideas used in the definition of volume and the proof of theorem 1.1.

By a \it 2-dimensional  angle structure \rm on a triangle, we mean an assignment of a number in $(0, \pi)$, called \it angle, \rm to each vertex of the
triangle. Classical geometric triangles are examples of 2-dimensional angle structures.  It is known that a 2-dimensional angle structure is the same as a Moebius triangle, i.e., a triangle of inner
angles in $(0, \pi)$ in the Riemann
sphere bounded by circles and lines (see [Lu2]).  We will interchange the use of terminology Moebius triangle and 2-dimensional 
angle structure in the rest of the paper.
  For a classical geometric tetrahedron,
the codimension-1 faces of them are classical geometric triangles of the same type. Similarly, for an angle structure on a
tetrahedron, the codimension-1 face of it is a 2-dimensional angle structure. Namely, the inner angle of a codimension-1
face is the spherical edge length  in the spherical vertex triangle. The edge length can be calculated by the Cosine Law for the spherical vertex triangle. To define the Schlaefli 1-form, we have 
to define the edge length 
of a Moebius triangle. This is achieved by  generalizing the Cosine Law for classical geometric triangles. The main observation is
that the side of the Cosine Law involving inner angles is still valid for 2-dimensional angle structures. We define the length of an edge in a 
Moebius triangle by using the Cosine Law. Having defined the lengths of edges of a Moebius
triangles, we are able to define the edge length of an angle structure on a tetrahedron by verifying the compatibility condition.
Namely, although each edge of a tetrahedron is adjacent to two triangular faces,  the length of the edge is independent 
of the choice two Moebius triangles adjacent it.   From this length, one forms the standard schlaefli 1-form defined on the space  $AS(3)$
of all angle structures and shows that the 1-form is closed. We define the volume to be the integration of this 1-form.

The volume function defined in this way automatically satisfies the Schlaefli formula. By the Schlaefli formula, we are able to 
identify the critical points of the volume in the space
of all angle structures $AS(M, T)$.  By the Lagrangian multiplier, a critical point $p \in AS(M, T)$ of the 
volume is the same as the following. If $e$ is an edge in $T$ adjacent to two tetrahedra $A$ and $B$, then the lengths of $e$ in
$A$ and $B$ (in their angle structures in $p$) are the same. This produces a condition for gluing angle structures on 3-simplexes
 along their codimension-1 faces. One the other hand, both 2-dimensional and 3-dimensional angle structures can be 
classified into three types: Euclidean, hyperbolic and spherical. The type of an angle structure depends only on 
the length of an edge. As a consequence, at the critical point, all 3-simplices have the same type.
 If each 3-simplex in $p$ is
a classical geometric tetrahedron, then condition (a) in theorem 1.1 holds. If there is a 3-simplex in $p$ which
is not a classical geometric tetrahedron, we are able to produce a normal surface so that it cuts
each 3-simplex in at most one normal disk and its Euler characteristic is
positive. The positivity of the Euler characteristic is a consequence of the following observation.
 Namely, the sum of the dihedral angles at two pairs of opposite edges of
a Euclidean or a hyperbolic tetrahedron is less than $2\pi$. Finally, if the triangulation has only one vertex, then
all edges become loops. Since the normal surface intersects some edge transversely in one point, it is non-separating.
\medskip
\noindent
1.4. The above finite dimensional variational set up can be extended easily to
find constant curvature cone metrics on 3-manifolds.

\medskip
\noindent
1.5.  
This paper is organized as follows.
In \S2, we establish a calculus Cosine Law and its derivative form. In \S3, we establish a generalized Schlaefli formula.
 In \S4, we define the edge lengths of 2-dimensional angle structures and state some of its consequences.
In \S5, we define the volume of angle structures on tetrahedra.
In \S6, we give a  classification of angle structures on tetrahedra. In \S7, we  sketch the proof of theorem 1.1.
 
\medskip
\noindent
1.6.  Acknowledgment.  We thank X.S. Lin,  Ben Chow,
 and Thomas Au for discussions on the related subjects. Bruno Martelli made a nice suggestion on improving the last statement in theorem 1.1.  Part of the work was done while  I was visiting
the Chinese University of Hong Kong and Yale University. We thank Thomas Au and Yair Minsky for their hospitality. 

\medskip
\noindent
{\bf \S2. A Calculus Cosine Law}
\medskip
\noindent
Recall that if a spherical triangle has three inner angles $x_1, x_2, x_3$ and edge lengths $y_1, y_2, y_3$ where $y_i$ is the
length of the edge opposite to the angle $x_i$, then the Cosine Law says,
$$ \cos(y_i) = \frac{ \cos x_i + \cos x_j \cos x_k}{\sin x_j \sin x_k} \tag 2.1$$
where $\{i,j,k\} =\{1,2,3\}$.  It turns out  (2.1) encodes all Cosine Laws for classical geometric triangles, i.e.,
 hyperbolic, Euclidean triangles and even hyperbolic right-angled hexagons if we interprete the terms appropriately. 
 To be more precise,
if a classical geometric triangle  in $K^2$ ($K^2 = S^2, E^2$ or $H^2$) has inner
angles $x_1, x_2, x_3$, then the length $y_i$ of the edge opposite to the
angle $x_i$ satisfies the Cosine Law,
$$ \cos (\sqrt{\lambda} y_i) = \frac{\cos x_i + \cos x_j \cos x_k}{\sin (x_j) \sin (x_k)} \tag 2.2$$
where $\{i,j,k\} =\{1,2,3\}$ and $\lambda =1,0,-1$ is the curvature of the space
$K^2$. 
 This prompts us to look at (2.1) from  analysis point of view.

\medskip
\noindent
2.1.
Suppose we have a single valued complex analytic function $y= y(x)$ where $y=(y_1, y_2, y_3) \in \bold C^3$ and 
$x=(x_1, x_2, x_3)$ is in some  open  connected set $\Omega$ in $\bold C^3$ so that $y$ and $x$ are
related by  (2.1). 

\medskip
\noindent
{\bf Theorem 2.1.} \it Suppose $\Omega \subset \bold C^6$ contains a diagonal point $(a,a,a)$ so that $y(a,a,a) = (b,b,b)$.
Let the indices $\{i,j,k\}$ be $\{1,2,3\}$ and  $A_{ijk} = \sin y_i \sin x_j \sin x_k$. Then

 $$ \quad \quad \quad  A_{ijk} = A_{jki}.  \tag a  $$

$$A^2_{ijk} =  1- \cos^2 x_i -\cos^2 x_j -\cos^2 x_k - 2 \cos x_i \cos x_j \cos x_k.  \tag b$$
At the point $x$ where $A_{ijk} \neq 0$, the following three hold, 
$$  \partial y_i/\partial x_i = \sin x_i /A_{ijk},    \tag c$$

$$   \partial y_i/\partial x_j = \partial y_i/\partial x_i \cos y_k, \tag d $$ 

$$  \cos(x_i) = \frac{ \cos y_i - \cos y_j \cos y_k}{\sin y_j \sin y_k}. \tag e$$
\rm

\noindent
2.2. {\bf Remarks. 1.} Formula (a)  is  sometimes called the \it Sine Law. \rm 

\noindent
{\bf 2.}  If we make a change of variable that
$z_i=\pi-y_i$, then (2.1) and (e) show that the map $z=z(x)$ is an involution, i.e., $z (z(x)) = x$. This corresponds to the
 duality theorem for spherical triangles. Using this, we see
that the partial derivatives $\partial x_i/\partial y_i$ and $\partial x_i/\partial y_j$ can be derived easily from theorem 2.1.
To be more precise, we have $\partial x_i/\partial y_j  = -\partial x_i/\partial y_i \cos x_k$.

\noindent
{\bf 3.} Let $F_1(x)= (x_1, \pi-x_2, \pi-x_3)$, then  $(\cos (y_1(F_1(x))), \cos (y_2(F_1(x))),$
 $\cos (y_3(F_1(x))))=(\cos(y_1(x)), \cos(\pi -y_2(x)), \cos(\pi-y_3(x))).$ This symmetry of the
Cosine Law (2.1) will be used extensively in the paper.

\medskip
\noindent
2.3. {\bf Proof.}
For simplicity, let $c_i = \cos x_i$ and $s_i = \sin x_i$ for $i=1,2,3$. Then by definition,
$$ A^2_{ijk} = \sin^2 y_i \sin^2 x_j \sin^2 x_k = (1-\cos^2 y_i) \sin ^2 x_j \sin^2 x_k$$
$$ = ( 1-c_j^2)(1- c_k^2) - (c_i + c_j c_k)^2$$
$$ = 1 - c_j^2 - c_k^2 + c_j^2 c_k^2 - c_i^2 - 2 c_i c_j c_k - c_j^2 c_k^2$$
$$ = 1- c_i^2 -c_j^2 -c_k^2 - 2c_i c_j c_k.$$
This shows that (b) holds. Now consider the analytic function $A_{ijk}/A_{jki}$. By (b), it takes values $\pm 1$.
By the assumption that $y(a,a,a)=(b,b,b)$, we see that value $1$ is achieved. Thus $A_{ijk}=A_{jki}$ in the
connected set $\Omega$. This shows that (a) holds.

By taking derivative of (2.1), we obtain,
$$ -\sin y_i \partial y_i/\partial x_i =  - \sin x_i /( \sin x_j \sin x_k).$$ 
This establishes part (c).

To see part (d), we take the derivative of (2.1) with respect to $x_j$. Thus
$$ -\sin y_i \partial y_i/\partial x_j =
 (1/s_k) [ (-s_j \cos x_k) s_j - \cos x_j ( \cos x_i + \cos x_j \cos x_k)]/ s^2_j$$
$$ = 1/(s^2_j s_k) (- \cos x_k - \cos x_i \cos x_j)$$
$$ =[ - s_i  /(s_j s_k)] [ \cos x_k + \cos x_i \cos x_j]/(s_i s_j)$$
$$= - (s_i /(s_j s_k)) \cos y_k$$
By dividing it by  $-\sin y_i$ and using part (c), we obtain identity (d).

\medskip
Finally, we derive (e).   By  (a) and (b), we have
$$ (\cos y_i - \cos y_j \cos y_k)/(\sin y_j \sin y_k)$$ $$ = [ (c_i + c_j c_k)/(s_j s_k) - ( c_j + c_i c_k)(c_k + c_i c_j) /( s_i s_k s_i s_j)]
/(\sin y_j \sin y_k))$$
$$=[ (1-c_i^2) ( c_i + c_j c_k) - (c_j c_k + c_i c_j^2 + c_i c_k^2 + c_i^2 c_j c_k)]/[(s_i s_k \sin y_j )(s_i s_j \sin y_k )]$$
$$ =( c_i + c_j c_k - c_i^3 - c_jc_kc_i^2 - c_j c_k - c_i c_j^2 - c_i c_k^2 - c_i^2 c_j c_k)/(A_{ikj}A_{ijk})$$
$$= c_i(  1- c_i^2 -c_j^2 -c_k^2 - 2c_i c_j c_k)              /A_{ijk}^2$$
$$=\cos x_i.$$

\medskip
\noindent
\S 3.  {\bf A Generalized Schlaefli Identity}

\medskip
\noindent
 In this section, we will fix a branch of the function $\arccos z$. Let $U$ be the set
 $\{ z \in \bold C | 0 < Re(z) < \pi\} \cup \{\sqrt{-1} x | x \in \bold R_{ \geq 0}\} \cup \{ \pi - \sqrt{-1}x | x \in \bold R_{\geq 0}
\}$. Then the restriction map $\cos z : U \to \bold C$ is a bijection. We define $\arccos z: \bold C \to U$ to
be the inverse and call it the principal branch of $\arccos z$. Note that $\arccos z$ is analytic in
 $\bold C -\{ x \in \bold R | |x| \geq 1\}$.  

\medskip
\noindent
3.1. Given a 3-simplex with vertices $\{v_1, ..., v_4\}$, a \it complex weight \rm on the 3-simplex is
an assignment of a complex number $x_{ij}=x_{ji}  \in \bold C -\{ \pi n | n \in \bold Z\}$ to the edge $v_iv_j$.
We consider $x_{ij}$ as a complex valued "dihedral angle".
 In this setting, the "face angle" $y^i_{jk}$ of the weighted 3-simplex in the face triangle $\Delta v_iv_jv_k$ at the vertex $v_i$
is defined to be the unique complex number $y^i_{jk}$ in $U$ so that

$$ \cos y^i_{jk} = \frac{ \cos x_{il}+ \cos x_{ij}\cos x_{ik} }{\sin x_{ij} \sin x_{ik}}  \tag 3.1$$
where $i,j,k,l$ are pairwise distinct.  We will assume that indices $i,j,k,l$ are always pairwise distinct in the sequel.
Note that if $\{x_{12}, ..., x_{34}\}$ forms the dihedral angles of a classical geometric tetrahedron,
then $y^i_{jk}$ is  the inner angle at $v_i$ in the triangle
$\Delta v_iv_jv_k$ by the Cosine Law (2.1).

\medskip
\noindent
{\bf Proposition 3.1.}(Compatibility)  \it Suppose the complex weight $x=(x_{12},...,x_{34})$  satisfies $\sin (y^i_{jk}) \neq 0$ for all indices. Then

$$ \frac{ \cos (y^i_{jk}) + \cos(y^j_{ik}) \cos (y^k_{ij})}{ \sin (y^j_{ik}) \sin (y^k_{ij})}
= \frac{ \cos (y^l_{jk}) + \cos(y^j_{lk}) \cos (y^k_{jl})}{ \sin (y^j_{lk}) \sin (y^k_{jl})}. \tag 3.2$$

\rm
The underlying geometric meaning of this proposition is that the length of the edge $v_j v_k$ in a classical geometric tetrahedron
can be calculated from any of the two triangles $\Delta v_i v_j v_k$ or $\Delta v_l v_j v_k$. 
The proof is a direct computation using the Sine Law in theorem 2.1. 

Note that the condition $\sin (y^i_{jk}) \neq 0$ for all indices $i,j,k$ is equivalent to $x_{ij} \pm x_{ik} \pm x_{il}
\neq (2n+1) \pi$ for some integer $n$ for all $\{i,j,k,l\}=\{1,2,3,4\}$. We call a complex weight $(x_{12},...,x_{34})$ \it non-degenerate \rm if  $\sin(y^i_{jk}) \neq 0$ for all indices.
In this case, 
the common value in (3.2) is denoted by $\cos z_{jk}$ where $z_{jk} \in U$. We call $z_{jk}$ the \it complex
length \rm of the edge $v_jv_k $. If the dihedral angles $(x_{12}, ..., x_{34})$ form
 the inner angles
of a spherical or a hyperbolic 3-simplex, by the Cosine Laws,
$\lambda \sqrt{\lambda} z_{ij}$ is equal to the length of edge $v_iv_j$ in the 3-simplex in $S^3$ or $H^3$ where $\lambda = \pm 1$ is the curvature
of $S^3$ or $H^3$.

\medskip
\noindent
3.2. 
For a non-degenerated complex weight $x=(x_{12}, ..., x_{34})$, the complex length $z_{ij} = z_{ij}(x)$ is a complex analytic
function of the weight $x$. 

\medskip
\noindent
{\bf Theorem 3.2.} \it For a non-degenerate complex weight $(x_{12},...,x_{34})$, then
  $$\partial z_{ij}/\partial x_{ab} = \partial z_{ab}/\partial x_{ij},$$
where  the indices  satisfy  $a \neq b$ and $i \neq j$.
In particular, the differential 1-form $\sum_{i,j} z_{ij} dx_{ij}$ is closed.
\rm
\medskip
The proof is computational using theorem 2.1.

\medskip
\noindent
\S4. {\bf  Lengths of Edges in  Moebius Triangles}

\medskip
\noindent
 Given a Moebius triangle of inner angles $(x_1, x_2, x_3) \in (0, \pi)^3$, we will define the
Moebius length (or length for simplicity) of an edge in the triangle in this section. 
The definition of the length is guided by formula (2.2). 

\medskip
\noindent
4.1. As a convention in this section, the function $\cos z$ is considered as a homeomorphism from 
the 1-dimensional subset $\bold L =\{ \sqrt{-1}x | x \in \bold R_{\geq 0}\}
\cup \{ x \in \bold R | x \in [0, \pi]\} \cup \{ \pi - \sqrt{-1} x
| x \in \bold R_{\geq 0}\}$ to the real line $\bold R$.  Let
$ \phi: \bold R \to \bold L $ be the homeomorphism 
given by $\phi(x) = x$ if $x \in [0, \pi]$; $\phi( x) = -\sqrt{-1} x$ if
$x \leq 0$; and $\phi(x) = \pi+\sqrt{-1}( \pi-x)$ if $x \geq \pi$.  In particular, the homeomorphism $f(x) =\cos (\phi(x)):
\bold R \to \bold R$ is given by $f(x) =  \cosh(x)$ when $x \leq 0$; $f(x) = \cos(x)$ when $ x \in [0, \pi]$; and
$f(x) = -\cosh(x-\pi)$ when $x \geq \pi$. The function $f$ is $C^1$-smooth and $f^{-1}$ is
continuous.
\medskip
\noindent
 {\bf Definition 4.1.}  Given a Moebius triangle of inner angles $(x_1, x_2, x_3)$, the \it length \rm (or more precisely
 the \it Moebius length) \rm  $z_i$ of the edge opposite
to the angle $x_i$ is the real number so that
$$ \cos ( \phi(z_i)) =\frac{ \cos x_i + \cos x_j \cos x_k}{ \sin x_j \sin x_k}  $$
where $\{i,j,k\} =\{1,2,3\}$.

\medskip

By formula (2.2), if $(x_1, x_2, x_3) \in (0, \pi)^3$ is a classical geometric triangle in $K^2 =S^2, H^2, E^2$ of curvature
$\lambda =1,0,-1$ so that lengths of the edges are $l_1, l_2, l_3$ measured in $K^2$ geometry, then the Moebius lengths of it
are $(\lambda l_1, \lambda l_2, \lambda l_3)$, i.e., the Moebius length is the signed length for classical geometric triangles.

\medskip
\noindent
4.2. Given a Moebius triangle $x =(x_1, x_2, x_3) \in (0, \pi)^3$, the i-th \it flip \rm $F_i(x)$ of $x$ is the Moebius triangle with inner angles
$y=(y_{a})$ where $y_i=x_i$ and $y_j=\pi-x_j$ for $j \neq i$.  

Let $E(2), H(2), S(2)$ be the subspaces of $AS(2)$ corresponding the Euclidean, hyperbolic and
spherical triangles respectively. Then evidently $F_i:S(2) \to S(2)$. It can be shown that for any $x \in AS(2)$, either
$x \in E(2) \cup H(2) \cup S(2)$ or there is a flip $F_i(x)$ so that $F_i(x) \in E(2) \cup H(2)$. We classify the Moebius
triangles into  Euclidean type, hyperbolic type or spherical type according to $x$ or $F_i(x)$ is the classical geometric triangle of the same
type. This classification is invariant under the flip operation.  
The following result summerizes the basic properties of the Moebius lengths, the type of Moebius triangles and flip operations.

\medskip
\noindent
{\bf Proposition 4.2.} \it The i-th edge length function $z_i: AS(2)=(0,\pi)^3  \to \bold R^1$ is a continuous function.
If $x$ is a classical geometric triangle in $K^2$ of curvature $\lambda=-1,1,0$, then the Moebius length $z_i(x)$ is $\lambda l_i$ where
$l_i$ is the length calculated in the classical geometry. Furthermore,

(a) For a non-Euclidean type Moebius triangle, the three Moebius lengths determine the three inner angles.

(b) If $F_i(x)$ is the i-th flip of $x$, then the Moebius length and the flip are related by
$$ z_j(F_i(x))= \pi - z_j(x)  \tag 4.1$$ 
$$ z_i(F_i(x)) = z_i(x)   \tag 4.2$$

(c) The type of a Moebius triangle is determined by one edge length. Namely, $x$ is of Euclidean type if and only
 if $z_i(x) \in \{0, \pi \}$;
$x$ is of spherical type if and only if $z_i(x) \in (0, \pi)$; and $x$ is of hyperbolic type if and only if $z_i(x) \in (-\infty, 0) \cup (\pi, \infty)$.

(d) A Moebius triangle is not a classical geometric triangle if and only if one edge  length is at least $\pi$.
\rm

\medskip
Part (a) follows from theorem 2.1 (e). Part (b) is a consequence of remark 2.2.3.
Parts (c) and (d) follow
from (b) and the definition.  There are exactly two edges in a non-classical geometric triangle of length at least $\pi$.
The flip operation about the vertex which is the intersection of the two edges of length $\geq  \pi$ is a classical geometric
triangle.

\medskip
\noindent
4.3. {\bf Remark.} One can take identities (4.1) and (4.2)  as the definition of the Moebius length of edges in Moebius triangles.

\medskip
\noindent
\S5. {\bf Volume of an Angle Structure on a Tetrahedron }
\medskip
\noindent
5.1. In this section, we assume that the indices $\{i,j,k,l\}=\{1,2,3,4\}$. An angle structure $x=(x_{12}, ..., x_{34})$
on a 3-simplex 
with vertices $v_1, ...,v_4$ has dihedral
 angle $x_{ij} = x_{ji}$ at the edge $v_iv_j$.  The space of all angle structures on a tetrahedron 
 $AS(3)$ is 
 $ \{ x \in(0, \pi)^6 | x_{ij}+x_{ik}+x_{il} > \pi, x_{ij}+x_{ik}-x_{il} < \pi, \text{ for all  i,j,k,l} \}$. 
 Given $x \in   AS(3)$ and a codimension-1 face $\Delta v_iv_jv_k$ of the tetrahedron, the \it induced 2-dimensional angle structure \rm
on $\Delta v_iv_jv_k$ is obtained by  assigning the positive number $y^i_{jk} \in (0, \pi)$ to the vertex $v_i$ where
$y^i_{jk}$ satisfies (3.1). Geometrically, if we construct a spherical triangle of inner angles
$x_{ij}, x_{ik}, x_{il}$ associated to the vertex $v_i$, then $y^i_{jk}$ is the spherical length of the edge opposite to
 the angle $x_{il}$ in the
spherical triangle.  There are four Moebius triangles $\Delta v_i v_j v_k$ appeared as codimension-1 faces of the tetrahedron. 
These
Moebius triangles have  Moebius lengths at the edges.
 By the definition of Moebius lengths
and compatibility proposition 3.1, we have,

\medskip
\noindent
{\bf Lemma 5.1.} \it The Moebius lengths of the edge $v_i v_j$ in the Moebius
triangles $\Delta v_iv_jv_k$ and $\Delta v_i v_j v_l$ are the same. \rm

\medskip
We call the common value in lemma 5.1 the \it length \rm (or Moebius length)
of the edge $v_i v_j$ in $x$, and denote it by $l_{ij}$.  Note that if the angle structure is a classical geometric tetrahedron
in the space $K^3=S^3, H^3, E^3$ of curvature $\lambda=1,-1,0$, then the Moebius length is $\lambda l$ where $l$ is the
length measured in the classical geometry $K^3$.

\medskip
\noindent
5.2.
Note that the regular Euclidean tetrahedron has the dihedral angle $\arccos(1/3)$.  A continuous differential 1-form
$\omega$ defined on a smooth manifold is said to be closed if its integration along each piecewise smooth loop is zero.

\medskip
\noindent
{\bf Theorem 5.2.} \it The continuous differential 1-form $\omega=1/2\sum_{i>j} l_{ij} dx_{ij}$ on the open convex 
polytope $  AS(3)$ is closed. Its integration
$$V(x) =\int_{a}^x \omega $$ where
$a = \arccos(1/3)(1,1,..,1)$ is a $C^1$-smooth function on $  AS(3)$, called the
volume (or Moebius volume). The function $V(x)$ has the following properties,

\medskip
\noindent
(a) (Schlaefli formula) $$\partial V(x)/\partial x_{ij} = l_{ij}/2,  \tag 5.1$$
\medskip
\noindent
(b) if $x \in   AS(3)$ is a classical geometric 3-simplex in the space of constant curvature $\lambda$ ($\lambda =-1,0,1$),
 then the volume $V(x) = \lambda^2 vol(x)$ where
$vol$ is the volume measured in the classical geometry of constant curvature $\lambda$. 
\medskip
\noindent
(c) The volume function $V$ can be extended continuously to the compact closure
 $\bar{  AS(3)}=\{ x \in [0, \pi]^3 | x_{ij}+x_{ik}+x_{il} \geq \pi, x_{ij}+x_{ik}-x_{il}  \leq \pi\}$.

\rm

\medskip
\noindent
\noindent
\S6. {\bf The Classification of Angle Structures and Flip Operations }
\medskip
\noindent
For an angle structure on a tetrahedron $x =(x_{rs}) \in AS(3)$, the i-th \it flip \rm $F_i(x) = (y_{rs})$ is the angle structure so that
$y_{ij} = x_{ij}$ and $y_{jk} = \pi-x_{jk}$ where $\{i,j,k,l\}$ are pairwise distinct.
It follows from the
definition and remark 2.2.3 that  the codimension-1 faces $\Delta v_iv_jv_k$  of $F_i(x)$ and $x$ are related by
the i-th flip and the faces $\Delta v_jv_kv_l$ of $x$ and $F_i(x)$ are
the same.
The corresponding vertex links of $F_i(x)$ and
$x$, considered as spherical triangles,  are either the same or related by a flip.  

\medskip
\noindent
6.1. Let $E(3), H(3)$ and $S(3)$ be the subspaces of $AS(3)$ corresponding to the Euclidean, hyperbolic, and spherical tetrahedra
respectively. It can be shown that,
\medskip
\noindent
{\bf Proposition 6.1.} \it For any $x \in AS(3)$, exactly one of the following holds,

(a) $x$ is a classical geometric tetrahedron, i.e., $x \in E(3) \cup H(3) \cup S(3)$,

(b) there is a flip $F_i$ so that $F_i(x) \in E(3) \cup H(3)$,

(c) there are two flips $F_i, F_j$, $i \neq j$, so that $F_iF_j(x) \in E(3) \cup H(3)$. \rm

\medskip
If an angle structure $x$ is obtained from a single flip $F_i(y)$ where
$y \in E(3) \cup H(3)$, then the edge length $l_{ij}$ of $x$ is at least $\pi$ and all other edge lengths
are at most $0$ by (4.1) and (4.2).
If $x$ is obtained from a double-flip $F_iF_j(y)$ where $y \in E(3) \cup H(3)$ with $i \neq j$, then
the edge lengths and dihedral angles of $x$ satisfy: (1) $l_{ab} \geq \pi$, $x_{ab} =\pi-y_{ab}$
for $\{a,b\}=\{i,k\}, \{i,l\}, \{j,k\}, \{j,l\}$; and (2) $l_{ij}, l_{kl} \leq 0$ and
$x_{ij}=y_{ij}, x_{kl}=y_{kl}$.  

\medskip
\noindent
{\bf Lemma 6.2}. \it Suppose $x=F_iF_j(y)$ for $y \in E(3)\cup H(3)$ is a double flip
of a Euclidean or a hyperbolic tetrahedron $y$. Then the sum of the dihedral angles at the two pairs of opposite edges of length at least $\pi$ is greater than $2\pi$. \rm

\medskip
Indeed, by the classification above, the lemma is equivalent to the following statement about
Euclidean and hyperbolic tetrahedra.
Namely, if $a_1, a_2, a_3, a_4$ are dihedral angles at two pairs of opposite edges in a Euclidean or
hyperbolic tetrahedron, then $a_1 + a_2+a_3 +a_4 < 2\pi$. 
The proof of  it is a simple exercise in geometry.

\medskip
\noindent
6.2. {\bf Remark.}
The relationship between the flip operation and  the volume  is the following.
For any angle structure $x \in AS(3)$, 
$  V(F_i(x)) + V(x) = \pi/2( x_{ij} + x_{ik} + x_{il} - \pi) $. One can in fact use this as the
definition of the volume of angle structure. Using this formula, we prove the continuous
extension of the volume to the closure of $AS(3)$ using  [Lu1]. Rivin [Ri2] has now a simpler
proof of it and proved a more generalized result.

\medskip
\noindent
\S7. {\bf The Sketch of the Proof of Theorem 1.1.}

\medskip
\noindent
The proof goes as follows. 

\medskip
\noindent
7.1. Recall that the type of a Moebius triangle is  determined by one edge length by
proposition 4.2. Using proposition 3.1, it follows that the types of all Moebius triangles
appeared as faces of 
an angle structure on a tetrahedron are the same. We define the type of an angle structure on
a tetrahedron to be
the type of a codimension-1 face of it.  By definition, the type of an angle structure is determined
by the length of one edge. To be more precise,
if the length of the ij-th edge  $l_{ij}$ is in $\{0, \pi\}$, then it is of Euclidean type; if $l_{ij}$ is in $(0,\pi)$, then it is of
spherical type; and
if $l_{ij}$ is in $(-\infty, 0) \cup (\pi, \infty)$, then it is of hyperbolic type. 
\medskip
\noindent
7.2. If $p \in AS(M, T)$ is a critical point of the volume, then Schlaefli formula shows that the following holds. Namely, if
$e \in T^{(1)}$ is an edge in the triangulation so that $e$ is adjacent to two 3-simplexes $A, B$, then the Moebius lengths
of $e$ in $A$ and $B$ (in $p$) are the same:
$$ l_e(A,p) = l_e(B,p).  \tag 7.1$$
Indeed, suppose $a_1$ and $a_2$ are the dihedral angles of $p=(a_1, a_2,..., a_n)$ at the edge $e$ inside $A$ and $B$. Now consider the
deformation $r(t) = (t, -t, 0,...,0)+p$. Since $p=(a_1, a_2,..., a_n)$ is the critical point of $V$, it follows that $d/dt|_{t=0} V(r(t))=0$. But
by the Schlaefli formula (5.1), the derivative is $1/2( l_e(A, p) - l_e(B,p))$. Thus (7.1) follows. 

\medskip
\noindent
7.3. By combining 7.1 and 7.2, the types of all angle structures on 3-simplexes in $p$ are the same.

\medskip
\noindent
7.4. If one 3-simplex in $p$ is spherical, then all 3-simplexes in $p$ are spherical. They are all classical spherical tetrahedra so that
their faces can be glued isometrically by (7.1). By the definition of angle structure that the sum of dihedral angles is $2\pi$ at each edge,
 we obtain a spherical metric on the 3-manifold $M$.

\medskip
\noindent
7.5. If all 3-simplexes in $p$ are classical hyperbolic 3-simplexes, then the same argument as in 7.4 shows that the manifold $M$ has
a hyperbolic metric.

\medskip
\noindent
7.6. If all 3-simplexes in $p$ are classical Euclidean tetrahedra, then we claim that $p$ is not a local
maximum point. In fact, 
 the  critical point $p$ is a local minimum of the volume. This is due to the
following two facts. First, if $x$ is an angle structure on a tetrahedron sufficiently close to a Euclidean tetrahedron,
then $x$ is a classical geometric tetrahedron, i.e., $S(3) \cup E(3) \cup H(3)$ is open in $AS(3)$.  Second, by theorem 5.2(b), the volume
of spherical and hyperbolic tetrahedra are positive and the volume of a Euclidean tetrahedron is zero.
Thus, the volume of points in $AS(M,T)$ sufficiently close to $p$ are none negative. On the other
hand, the volume of $p$ is zero. This shows that the point $p$ is a local minimum. It is easy to
show that the volume function is not a constant on any open subset of $AS(M, T)$. Thus, the
critical point $p$ is not a local maximum. 

\medskip
\noindent
7.7. We claim,
\medskip
\noindent
{\bf Proposition 7.1.} \it If one 3-simplex in $p$ is not a classical geometric tetrahedron, then the triangulation $T$ contains
a normal surface $S$ of positive Euler characteristic which cuts each 3-simplex in at most one normal disk. \rm

\medskip
To prove this, let $X$ be the set of all edges in the triangulation so that its length is at least $\pi$ (in $p$).  This set is non-empty since
there are non-classical geometric tetrahedra in $p$. By definition, the intersection of $X$ with each 3-simplex $A$ in $T$ consists of
the following three cases:

 (1) $X \cap A =\emptyset$ if $A$ is a classical geometric tetrahedron in $p$,

(2) $X \cap A$ consists of three edges from a vertex if $A$ is $F_i(x)$ for a classical geometric tetrahedron $x \in E(3) \cup H(3)$,

(3) $X \cap A$ consists of pair of opposite edges if $A$ is $F_iF_j(x)$ for a classical geometric tetrahedron $x \in E(3) \cup H(3)$.

For each 3-simplex  $A$ in cases (2) and (3), we construct a normal disk in $A$ whose vertices are in $X$. These normal disks
form a normal surface $S$ in $T$. Note
that the normal surface $S$ intersects the each tetrahedron in $T$ in at most one normal disk. 
We claim the normal surface $S$ has positive Euler characteristic. To see this, let us consider the
CW-decomposition of the surface $S$ formed by the normal disks in $S$. We assign each 
vertex of each 2-cell in the CW-decomposition a number, called the \it inner angle\rm, which  is the corresponding dihedral angle in $p$. By the construction, these 
inner angles satisfy the following three conditions.
First, the 
sum of all inner angles at each vertex in $S$ is $2\pi$. Second, the sum of
the inner angles in each normal triangle greater than $\pi$. Finally, the sum of
all inner angles in each normal quadrilateral is greater than $2\pi$ by lemma 6.2.
By the Gauss-Bonnet theorem,  it follows that the normal surface $S$ has
positive Euler characteristic. 

\medskip
\noindent
\centerline{\bf References}
\medskip

[Gu] Gueritaud,  Francois, On canonical triangulations of the mapping tori over the punctured torus,
http://front.math.ucdavis.edu/math.GT/0406242.

[GL] Guo, Ren and Luo, Feng, On a conjecture of Minlor on the volume of simplexes, II, 
http://front.math.ucdavis.edu/math.GT/0510666

[JR]  Jaco, William; Rubinstein, J. Hyam, $0$-efficient triangulations of 3-manifolds. J. Differential Geom.
   65 (2003), no. 1, 61--168. 

[La]  Lackenby, Marc, Word hyperbolic Dehn surgery. Invent. Math. 140 (2000), no. 2, 243--282.
 
[Lu1]  Luo, Feng, Continuity of the volume of simplices in classical geometry, to appear,  Commun. Contemp. Math.
http://front.math.ucdavis.edu/math.GT/0412208. 

[Lu2]  Luo, Feng,  Triangulations in Moebius geometry. Trans. Amer. Math. Soc. 337 (1993), no. 1, 181--193.

[Lu3] Luo, Feng, Volume and angle structures on 3-manifolds, II, in preparation.

[LT] Luo, Feng and Tillman, Stephen, Angle structures and  normal surfaces, preprint, 2005.

[Mi] Milnor, J, The Schlaefli differential equality. In Collected papers, vol. 1.
Publish or Perish, Inc., Houston, TX, 1994.

[Mu] Murakami, Jun, Generalized volume and geometric structure of 3-manifolds,
http://www.f.waseda.jp/murakami/papers/VSrev.pdf

[Ri1] Rivin, Igor, Combinatorial optimization in geometry. Adv. in Appl. Math. 31 (2003), no. 1, 242--271.

[Ri2] Rivin, Igor, Continuity of volumes -- on a generalization of a conjecture of J. W. Milnor,
http://front.math.ucdavis.edu/math.GT/0502543

\medskip
Department of Math., Rutgers University, New Brunswick, NJ 08854, USA

email: fluo\@math.rutgers.edu

\end